# APPROXIMATE THEORY-AIDED ROBUST EFFICIENT FACTORIAL FRACTIONS UNDER BASELINE PARAMETRIZATION

Rahul Mukerjee  and  S. Huda
Indian Institute of Management Calcutta     Department of Statistics and OR
Joka, Diamond Harbour Road                  Faculty of Science, Kuwait University
Kolkata 700 104, India                      P.O. Box-5969, Safat-13060, Kuwait

**Abstract**: With reference to a baseline parametrization, we explore highly efficient fractional factorial designs for inference on the main effects and, perhaps, some interactions. Our tools include approximate theory together with certain carefully devised discretization procedures. The robustness of these designs to possible model misspecification is investigated using a minimaxity approach. Examples are given to demonstrate that our technique works well even when the run size is quite small.

*Key words*: A-criterion, bias, binary design, directional derivative, discretization, minimaxity, model misspecification, multiplicative algorithm, nonorthogonality, requirement set.

*Mathematics Subject Classification*: 62K15.

## 1. Introduction

Design of optimal or efficient factorial fractions has received significant attention in recent years; see Mukerjee and Wu (2006), Wu and Hamada (2009) and Xu, Phoa and Wong (2009) for surveys and further references. A vast majority of the existing work on fractional factorial designs centers around the well-known orthogonal parametrization (Gupta and Mukerjee, 1989) where orthogonal arrays play a key role in the construction of optimal fractions.

In recent years, however, a baseline parametrization has started gaining popularity in factorial experiments. It has found use in microarray experiments (Yang and Speed, 2002) and can also arise naturally in many other situations whenever each factor has a control or baseline level. Following Kerr (2006), an example is given by a toxicological study with binary factors, each representing the presence or absence of a toxin, where the state of absence is a natural baseline level of each factor. Also, in agricultural or industrial experiments, the currently used level of each factor can well constitute the baseline level. Under baseline parametrization, optimal paired comparison designs for full factorials have been studied by several authors in the context of microarray experiments; see Banerjee and Mukerjee (2008), Zhang and Mukerjee (2013), and the references therein. As for designing efficient or optimal factorial fractions under this parametrization, not much work has been reported so far beyond the results in Mukerjee and Tang (2012) and Li, Miller and Tang (2014) on the two-level case.

Indeed, the lack of orthogonality in baseline parametrization makes the combinatorics for general fractional factorials extremely involved – for example, unlike what happens in orthogonal parametriza-tion, beyond two-level factorials, orthogonal arrays cease to remain optimal for main effects even in the absence of interactions; see Mukerjee and Tang (2012). Given this difficulty, one naturally wonders if, pending the development of perfect optimality results, it is possible to have designs with assured high efficiency, say over 0.90 or 0.95, which would suffice in most applications. From this perspective, we develop a technique for finding, under the baseline parametrization, designs which are highly efficient under a given model and, at the same time, enjoy robustness to model misspecification. This is done via the approximate theory which in a sense exploits the underlying nonorthogonality to our advantage and provides a benchmark for assuring the efficiency of exact designs. It is, however, seen that discretization of the resulting optimal design measure via a commonly used rounding off technique can have disastrous consequences. As a remedy, we propose certain procedures that yield exact designs with high efficiency. We add that our approach to robustness attempts to address an issue raised by Yin and Zhou (2014) who worked with the orthogonal parametrization and stressed on the need for similar studies under nonorthogonal parametrizations, like baseline.

Section 2 introduces the baseline parametrization and formulates the design objective with reference to a model which is kept quite flexible. It includes the baseline effect, all main effects and perhaps some interactions. Then in Section 3, we show how approximate theory, in conjunction with certain discretization procedures, can yield highly efficient fractions of general factorials with an arbitrary number of levels for each factor and arbitrary run size. The robustness of these designs to model misspecification is explored in Section 4. Examples are given in Section 5 to demonstrate that our approach works well even when the run size is quite small.

**2. Baseline parametrization and design objective**

*2.1 Baseline parametrization*

Consider an $m_1 \times ... \times m_n$ factorial with $n$ factors $F_1,...,F_n$. For $1 \le i \le n$, the levels of $F_i$ are coded as $0, 1, ..., m_i - 1$, where 0 is the control or baseline level. Thus there are $v = m_1...m_n$ treatment combinations $j_1...j_n$ ($0 \le j_i \le m_i - 1, 1 \le i \le n$). Let $\tau(j_1...j_n)$ denote the treatment effect of the treatment combination $j_1...j_n$. Then the baseline parametrization for a full factorial model is given by

$$\tau(j_1...j_n) = \Sigma_{u_1 \in \{0, j_1\}} ... \Sigma_{u_n \in \{0, j_n\}} \theta(u_1...u_n), \qquad (2.1)$$

for every $j_1...j_n$, where the sum on each $u_i$ is over $u_i \in \{0, j_i\}$, i.e., $u_i$ is fixed at 0 if $j_i = 0$, while $u_i$ ranges over only 0 and $j_i$ if $j_i > 0$. In (2.1), $\theta(0...0)$ is the baseline effect, while $\theta(u_1...u_n)$ is a main or interaction effect parameter for every $u_1...u_n \ne 0...0$, depending on which $u_i$'s are nonzero. For instance, with a $3^2$ factorial, (2.1) yields

$$\tau(00) = \theta(00), \quad \tau(j_1 0) = \theta(00) + \theta(j_1 0) \ (j_1 = 1,2), \quad \tau(0 j_2) = \theta(00) + \theta(0 j_2) \ (j_2 = 1,2),$$

$$\tau(j_1 j_2) = \theta(00) + \theta(j_1 0) + \theta(0 j_2) + \theta(j_1 j_2) \ (j_1, j_2 = 1,2), \tag{2.2}$$

where $\theta(00)$ is the baseline effect, main effect $F_1$ is represented by $\theta(10), \theta(20)$, main effect $F_2$ is represented by $\theta(01), \theta(02)$, and interaction $F_1 F_2$ is represented by $\theta(11), \theta(12), \theta(21), \theta(22)$. It is easy to see that the parametrization (2.1) is nonorthogonal – if the $v$ equations (2.1) arising from the $v$ treatment combinations are written in matrix notation, then the columns of the coefficient matrix in the right-hand side are not mutually orthogonal even when they correspond to different factorial effects.

We begin by considering a reduced version of (2.1) which includes the baseline effect and a collection of possibly significant factorial effects, such as all main effects and perhaps some interactions as well. The factorial effects kept in the model constitute the requirement set, denoted by $R$. Write $\theta_0$ for the baseline effect $\theta(0...0)$, and $\theta$ for the vector of parameters representing the factorial effects in $R$. Interest lies in inference on $\theta$, while $\theta_0$ is a nuisance parameter. Let $q$ be the dimension of $\theta$. Thus in the $3^2$ example, if $R$ consists of main effects $F_1$ and $F_2$, then parameters representing $F_1 F_2$ do not appear in the model. So, (2.2) reduces to $\tau(j_1 j_2) = \theta(00) + \theta(j_1 0) + \theta(0 j_2)$ for $j_1, j_2 = 1, 2$, but remains unchanged when $j_1 = 0$ or $j_2 = 0$. In this case, $\theta = (\theta(10), \theta(20), \theta(01), \theta(02))^T$ consists of the main effects parameters, with the superscript $T$ denoting transposition, and $q = 4$.

In order to streamline the presentation, we now label the $v$ treatment combinations $1,...,v$ in the lexicographic order, e.g., in a $3^2$ factorial, the lexicographic ordering of the $v = 9$ treatment combinations is 00, 01, 02, 10, 11, 12, 20, 21, 22, and these are labeled 1, 2,..., 9, respectively. More generally, in an $m_1 \times ... \times m_n$ factorial, any treatment combination $j_1 ... j_n$ is labeled $k$ where

$$k = \mu_1 j_1 + ... + \mu_n j_n + 1, \tag{2.3}$$

and $\mu_i = v/(m_1...m_i)$, $1 \le i \le n$. With $k$ as above, we simply write $\tau_k$ for $\tau(j_1...j_n)$. Clearly, then for the model introduced above, $\tau_k = \theta_0 + z_k^T \theta$, where $z_k$ is a known $q \times 1$ vector, $1 \le k \le v$; e.g., in the $3^2$ example, the treatment combinations 10 and 21 are labeled 4 and 8, respectively, and if $R$ consists of the two main effects then, with $\theta$ as shown in the previous paragraph, $z_4 = (1, 0, 0, 0)^T$ and $z_8 = (0, 1, 1, 0)^T$. Let $\tau = (\tau_1,...,\tau_v)^T$. Then our model can be expressed as

$$\tau = \theta_0 1_v + Z\theta, \tag{2.4}$$

where $1_v$ is the $v \times 1$ vector of ones, and $Z$ is the $v \times q$ matrix with rows $z_1^T,...,z_v^T$.

In what follows, for any vector $b = (b_1,...,b_l)^T$, with $B = b_1 + ... + b_l > 0$, we define

$$D(b) = \text{diag}(b_1,...,b_l), \qquad \Delta(b) = \text{diag}(b_1,...,b_l) - B^{-1}bb^T. \qquad (2.5)$$

Also, for any positive integer $l$, we write $0_l$ for the $l \times 1$ vector of zeros, $1_l$ for the $l \times 1$ vector of ones, and $I_l$ for the identity matrix of order $l$.

*2.2 Design objective*

Consider now an $N$-run design $d$, consisting, say, of treatment combinations labeled $k(1),...,k(N)$ which are not necessarily distinct. Let $Y_d$ denote the $N \times 1$ observational vector arising from $d$. Assume that the observational errors are uncorrelated with a common variance $\sigma^2$. Then, corresponding to (2.4), we have the model

$$E(Y_d) = \theta_0 1_N + Z_d \theta, \qquad \text{cov}(Y_d) = \sigma^2 I_N, \qquad (2.6)$$

where $Z_d$ is the $N \times q$ matrix with rows $z_{k(1)}^T,..., z_{k(N)}^T$. Let $L_N = I_N - N^{-1} 1_N 1_N^T$. Then under (2.6), the information matrix for the parametric vector $\theta$ of interest is given by

$$H_d = Z_d^T L_N Z_d. \qquad (2.7)$$

We consider only designs $d$ which keep $\theta$ estimable, i.e., have a nonsingular $H_d$; one can check that the condition $N \geq q+1$ is necessary for this purpose. By (2.6), for such a design $d$, we have

$$\hat{\theta} = H_d^{-1} Z_d^T L_N Y_d, \qquad \text{cov}_d(\hat{\theta}) = \sigma^2 H_d^{-1}, \qquad (2.8)$$

where $\hat{\theta}$ is the best linear unbiased estimator of $\theta$.

We aim at finding a design $d$ so as to minimize the average variance of the elements of $\hat{\theta}$ among all $N$-run designs (the issue of model robustness is taken up later). By (2.8), this calls for minimizing $\text{tr}(H_d^{-1})$. So, we go by the *A*-criterion which is justified here because interest lies in the parameters in $\theta$ as they stand, and not in their linear functions; cf. Kerr (2012). As hinted in the introduction, this exact design problem is combinatorially very complex, due to the nonorthogonality of the baseline parametrization. Use of approximate theory, however, allows us to make considerable progress.

## 3. Approximate theory and discretization procedures

*3.1 Approximate theory*

To motivate the ideas, we consider an alternative expression for $H_d$. For $1 \leq k \leq v$, let the treatment combination labeled $k$ appear $r_{dk} (\geq 0)$ times in $d$, where $r_{d1} + ... + r_{dv} = N$. The design is called binary if each $r_{dk}$ is 0 or 1. Let $r_d = (r_{d1},...,r_{dv})^T$. Then from (2.5), (2.7) and the definition of $Z_d$,

$$H_d = Z_d^T L_N Z_d = Z^T \Delta(r_d) Z = NM(p^{(d)}), \qquad (3.1)$$

where $p^{(d)} = N^{-1} r_d$ and $M(p^{(d)}) = Z^T \Delta(p^{(d)}) Z$. The discreteness of $r_d$ induces the same on $p^{(d)}$, but in order to employ the approximate theory we, for now, treat the elements of $p^{(d)}$ as nonnegative continuous variables which add up to 1. Any such $p = (p_1, ..., p_v)^T$ is called a design measure assigning masses $p_1, ..., p_v$ on the treatment combinations $1, ..., v$. By (3.1), then the design problem reduces to finding an optimal design measure which minimizes $\phi(p)$ over all possible $p$, where

$$\phi(p) = \text{tr}\{M(p)\}^{-1}, \text{ if } M(p) \text{ is nonsingular,}$$
$$= +\infty, \qquad \text{otherwise,} \qquad (3.2)$$

with

$$M(p) = Z^T \Delta(p) Z. \qquad (3.3)$$

Due to elimination of the nuisance parameter $\theta_0$, the expression (3.3) for $M(p)$ is not linear in the elements of $p$. Lemmas 1 and 2 below show that the still the basic ideas of approximate theory continue to hold in our setup. The first of these establishes the convexity of $\phi(p)$ and the second one characterizes $p$ so as to minimize $\phi(p)$. Their proofs appear in the appendix.

**Lemma 1**: *For any design measures $p$ and $\tilde{p}$ and any $0 < \varepsilon < 1$,*

$$(1 - \varepsilon)\phi(p) + \varepsilon\phi(\tilde{p}) \geq \phi((1 - \varepsilon)p + \varepsilon\tilde{p}).$$

**Lemma 2**. *A design measure $p$ minimizes $\phi(p)$ if and only if $M(p)$ is nonsingular and*

$$(z_k - Z^T p)^T \{M(p)\}^{-1} \{M(p)\}^{-1} (z_k - Z^T p) \leq \text{tr}\{M(p)\}^{-1}, \ 1 \leq k \leq v.$$

Lemma 2 leads to a multiplicative algorithm for numerical determination of the optimal design measure. It looks like that in Zhang and Mukerjee (2013) who considered paired comparison designs for full factorials, but is more elaborate than theirs due to the nonlinearity of our $M(p)$. The algorithm starts with the uniform measure $p^{[0]} = (1/v, ..., 1/v)^T$, and finds $p^{[h]} = (p_1^{[h]}, ..., p_v^{[h]})^T$, $h = 1, 2, ...$ recursively as

$$p_k^{[h]} = p_k^{[h-1]} \frac{(z_k - Z^T p^{[h-1]})^T \{M(p^{[h-1]})\}^{-1} \{M(p^{[h-1]})\}^{-1} (z_k - Z^T p^{[h-1]})}{\text{tr}\{M(p^{[h-1]})\}^{-1}}, \quad 1 \leq k \leq v, \quad (3.4)$$

till a design measure $p^{[h]}$, satisfying

$$(z_k - Z^T p^{[h]})^T \{M(p^{[h]})\}^{-1} \{M(p^{[h]})\}^{-1} (z_k - Z^T p^{[h]}) - \text{tr}\{M(p^{[h]})\}^{-1} \leq t, \quad 1 \leq k \leq v, \quad (3.5)$$

is obtained, where $t$ ($> 0$) is a preassigned small quantity. If we denote the terminal design measure meeting (3.5) by $\hat{p} = (\hat{p}_1,...,\hat{p}_v)^T$, then arguments similar to those in Silvey (1980, p. 35) and the proof of our Lemma 2 show that $\phi(\hat{p}) \leq \bar{\phi} + t$, where $\bar{\phi}$ is the minimum of $\phi(p)$ over all possible $p$. We take $t = 10^{-10}$. Then $\hat{p}$ represents the optimal design measure with accuracy up to nine places of decimals, Even at this level of accuracy, the algorithm in (3.4) and (3.5) works quite fast. In all the examples in Section 5 and many others, it terminates almost instantaneously.

In view of the foregoing discussion, for an $N$-run exact design $d$, it follows from (3.1) and (3.2) that

$$\text{tr}(H_d^{-1}) \geq s/N, \tag{3.6}$$

where $s = \text{tr}\{M(\hat{p})\}^{-1} - t$, with $t = 10^{-10}$. Thus a lower bound to the efficiency of $d$ under model (2.6) is obtained as

$$\text{eff}_{lb} = s/\{N\text{tr}(H_d^{-1})\} \tag{3.7}$$

Of course, $t$ is so small that its inclusion or otherwise in $s$ has effectively no impact on (3.7), as seen in our examples and many others. We remark that the bound (3.7), relative to the optimal design measure, is typically unattainable in an exact setup. So, it acts as a conservative benchmark and the actual efficiency of $d$ among $N$-run exact designs is often higher than this bound.

*3.2 Discretization procedures*

Approximate theory exploits to our advantage the two features of the baseline parametrization that hinder the development of exact optimality results, namely, nonorthogonality and lack of symmetry between the baseline and other levels of the factors. Because of these features, much unlike what often happens under orthogonal parametrization, our optimal design measure $\hat{p}$ turns out to be far from uniform. For instance, in a $2^5 \times 3$ factorial, with the requirement set consisting of all main effects and interactions $F_1 F_6$, $F_2 F_6$, it assigns masses 0.0479, 0.0150 and 0.0057 to treatment combinations 000000, 110000 and 110011, respectively. This lack of uniformity provides useful guidance to finding good exact designs – e.g., it is intuitively clear that an efficient exact design should include those treatment combinations where the optimal design measure assigns greater masses.

Notwithstanding the above, serious difficulties remain in translating the optimal design measure $\hat{p} = (\hat{p}_1,...,\hat{p}_v)^T$ to efficient exact designs. Consider, for instance, the common rounding off technique where, with a given run size $N$, one (i) multiplies $\hat{p}_1,...,\hat{p}_v$ by a positive constant $c$ so chosen that the quantities $c\hat{p}_1,...,c\hat{p}_v$, when rounded off to nearest integers, say $r_1,...,r_v$, add up to $N$, (ii) obtains an exact design $d$ where the treatment combination labeled $k$ appears $r_k$ times, $1 \leq k \leq v$, and (iii) hopes $d$ to be highly efficient with the belief that $r_1,...,r_v$ are approximately proportional to $\hat{p}_1,...,\hat{p}_v$. This often fails to work, especially for smaller run sizes which are of practical interest from the viewpoint of experimental economy. First, there may not exist any constant $c$ such that the resulting $r_1,...,r_v$ add up to $N$. Second, even when such a $c$ exists, the design $d$ obtained as above may turn out to be of poor efficiency, if not outright singular. Thus in the $2^5 \times 3$ example mentioned above, it is seen that rounding off fails to produce a nonsingular exact design for any $N \leq 32$.

For large $N$, on the other hand, the consequences of rounding off are negligible compared to $N$, and hence rounding off works well. This enables us to devise procedures for discretization of $\hat{p}$ in such a manner that the resulting exact designs retain high efficiency even for smaller $N$.

**Procedure A**:

I. Start with an $N_1$-run exact design $d(N_1)$, obtained from $\hat{p}$ by rounding off, such that $d(N_1)$ has very high efficiency, say with $\text{eff}_{lb} \geq 0.98$, where $\text{eff}_{lb}$ is given by (3.7). Here $N_1$ can be much larger than the target run size $N$ and $d(N_1)$ can be non-binary.

II. For $i = 0,1,..., N_1 - N - 1$, obtain an $(N_1 - i - 1)$-run design $d(N_1 - i - 1)$ from an $(N_1 - i)$-run design $d(N_1 - i)$ as follows:

(a) Consider all possible deletions of one run from $d(N_1 - i)$. Among the resulting $(N_1 - i - 1)$-run designs, let $d^*$ be the one with largest $\text{eff}_{lb}$. If $\text{eff}_{lb} \geq 0.95$ for $d^*$, then take $d(N_1 - i - 1)$ as $d^*$.

(b) Else, consider all possible deletions of two runs from $d(N_1 - i)$ coupled with all possible additions of one run from amongst the $v$ treatment combinations. Take $d(N_1 - i - 1)$ as the design with the largest $\text{eff}_{lb}$ among all $(N_1 - i - 1)$-run designs so generated.

The above procedure attempts to guard against $\text{eff}_{lb}$ going below 0.95 at any stage. Even if this is not always achieved, the final design $d(N)$, of run size $N$, typically has high efficiency, with $\text{eff}_{lb} \geq$ 0.90, even when $N$ is rather small. Several examples in Section 5 illustrate this point. Our computations suggest that Zhang and Mukerjee's (2013) approach, involving only deletion or only addition of runs, does not often work here in ensuring high efficiency for smaller $N$. This why, our procedure A is different from theirs, allowing both deletion and addition. Indeed, one may wonder about various modifications of A, such as deleting three runs and adding two in Step II (b), but we find that these only increase the computational burden significantly without entailing much gain in efficiency.

Since the initial design $d(N_1)$ in procedure A is often non-binary, the final output $d(N)$ can also be so, especially when $N$ is not too small. However, as seen in the next section, non-binary designs are unlikely to be robust to model misspecification. So, we also consider two variants of this procedure which always lead to binary designs while aiming at high efficiency.

**Procedure B1**:

I. Start with the full factorial design where each of the $v$ treatment combinations appears once. Set $N_1 = v$, and call this design $d(N_1)$.

II. Same as in A with the only change that if possibility (b) arises then the added run is so chosen that only binary designs are entertained. Note that starting with a binary $d(N_1 - i)$, one necessarily gets a binary $d(N_1 - i - 1)$ if (a) arises.

**Procedure B2**:

I. Start with the full factorial design. Set $N_1 = v$, and call this design $d(N_1)$.

II. For $i = 0, 1, \ldots, N_1 - N - 1$, consider all possible deletions of one run from $d(N_1 - i)$. Take $d(N_1 - i - 1)$ as the design with the largest $\text{eff}_{lb}$ among all $(N_1 - i - 1)$-run designs so generated.

While the optimal design measure $\hat{p}$ has no role in Step I of B1 and B2, it does influence the final outcome through the use of $\text{eff}_{lb}$ in Step II. Procedure B2 is very fast since it involves only deletion. Hence there is no harm in trying it first and if one is not happy with the $\text{eff}_{lb}$ value of the resulting $d(N)$, one can employ B1 which is not much time consuming either. In fact, B2 can help in making B1 even faster. One can employ B2 to quickly find an $N_1$-run binary design $d(N_1)$, with $\text{eff}_{lb} \geq 0.98$, and initiate B1 from this $d(N_1)$ rather than the full factorial. Typically, such $N_1$ is large, but still smaller than $v$, and this may expedite B1 in some situations. Illustrative examples follow in Section 5.

**4. Robustness to model misspecification**

With reference to the model (2.4), write $X = [1_v, Z]$. Let $C(X)$ denote the column space of $X$ and $C^\perp(X)$ the orthocomplement thereof in the $v$-dimensional Euclidean space. The dimensions of $C(X)$ and $C^\perp(X)$ are $q+1$ and $v-q-1$, respectively, since $X$ has full column rank, a fact which is not hard to see from a matrix representation of the baseline parametrization; cf. Banerjee and Mukerjee (2008). Model (2.4) amounts to assuming that $\tau \in C(X)$. In the absence of this model assumption, $\tau$ is any vector in the $v$-dimensional Euclidean space and hence can be represented as

$$\tau = \theta_0 1_v + Z\theta + P\xi, \tag{4.1}$$

where $P$ is a $v \times (v-q-1)$ matrix whose columns form an orthonormal basis of $C^{\perp}(X)$ and $\xi$ is a column vector of order $v-q-1$. The choice of $P$ is non-unique but, as seen below, this does not affect our findings.

Corresponding to (4.1), the model (2.6) for an $N$-run exact design $d$ gets modified to

$$E(Y_d) = \theta_0 1_N + Z_d \theta + P_d \xi, \qquad \text{cov}(Y_d) = \sigma^2 I_N, \qquad (4.2)$$

where $P_d$ is related to $P$ exactly in the same way as $Z_d$ is related to $Z$. Under (4.2), $\hat{\theta}$ in (2.8) no longer remains unbiased for $\theta$. It has now bias $H_d^{-1} Z_d^T L_N P_d \xi$ and its mean squared error matrix equals

$$\text{MSE}_d(\hat{\theta}) = \sigma^2 H_d^{-1} + H_d^{-1} Z_d^T L_N P_d \xi \xi^T P_d^T L_N Z_d H_d^{-1}. \qquad (4.3)$$

A model robust efficient design aims at keeping $\text{tr}\{\text{MSE}_d(\hat{\theta})\}$ small, a problem which is complicated by the fact that $\xi$ is unknown. To overcome this difficulty, in the spirit of Mukerjee and Tang (2012), we adopt a Bayesian inspired approach without a very explicit prior specification. Motivated by the hope that even if the assumed model (2.4) is inadequate, it is not too far from the true model (4.2), we consider the class $\Pi$ of priors $\pi$ such that

$$\lambda_{\max}\{E_\pi(\xi \xi^T)\} \leq \delta^2, \qquad (4.4)$$

where $\lambda_{\max}$ stands for largest eigenvalue and $E_\pi$ denotes expectation with respect to $\pi$. With a view to protecting against the worst scenario over all possible priors $\pi$ satisfying (4.4), we proceed to find a design $d$ which keeps

$$\psi_d = \max_{\pi \in \Pi} E_\pi[\text{tr}\{\text{MSE}_d(\hat{\theta})\}] \qquad (4.5)$$

small. This minimaxity approach is inspired by Wilmut and Zhou (2011), Lin and Zhou (2013) and Yin and Zhou (2014), who worked with the orthogonal parametrization. However, although (4.4) has essentially the same motivation as the condition that they imposed, it is not identical to theirs and seems to be more suited to the present context. We return to this point later but, at this stage, remark that (4.4) invariant of the choice of $P$. To see this, note that if $\tilde{P}$ is another matrix with orthonormal columns and the same column space as $P$, then $P\xi = \tilde{P}\tilde{\xi}$ for some $\tilde{\xi}$, and $P = \tilde{P}K$ for some orthogonal matrix $K$. Therefore, $\tilde{\xi} = K\xi$ and (4.4) is equivalent to its counterpart with $\xi$ changed to $\tilde{\xi}$. Lemma 3, proved in the appendix, helps in finding $\psi_d$ explicitly.

**Lemma 3**. (a) *Let* $V_d = H_d^{-1} Z^T \Delta(r_d) \Delta(r_d) Z H_d^{-1}$ *and* $W = \{Z^T \Delta(1_v) Z\}^{-1}$. *Then*

$$H_d^{-1} Z_d^T L_N P_d P_d^T L_N Z_d H_d^{-1} = V_d - W.$$

(b) *The matrix $V_d - H_d^{-1}$ is nonnegative definite (nnd), and $V_d = H_d^{-1}$ if d is binary.*

By (4.4), $\delta^2 I_{v-q-1} - E_\pi(\xi\xi^T)$ is nnd for any $\pi \in \Pi$. Hence by (4.3), for any such $\pi$,

$$E_\pi[\text{tr}\{\text{MSE}_d(\hat{\theta})\}] \leq \sigma^2 \text{tr}(H_d^{-1}) + \delta^2 \text{tr}\{H_d^{-1} Z_d^T L_N P_d P_d^T L_N Z_d H_d^{-1}\}, \qquad (4.6)$$

and this upper bound is attained if $E_\pi(\xi\xi^T) = \delta^2 I_{v-q-1}$ which also meets (4.4). Thus $\psi_d$ in (4.5) equals right-hand side of (4.6), and invoking Lemma 3(a), it follows that

$$\psi_d = \sigma^2 \text{tr}(H_d^{-1}) + \delta^2 \{\text{tr}(V_d) - \text{tr}(W)\}, \qquad (4.7)$$

which is invariant of the choice of *P*. By (3.6), (4.7) and Lemma 3(b),

$$\psi_d \geq (\sigma^2 + \delta^2)\text{tr}(H_d^{-1}) - \delta^2 \text{tr}(W) \geq (\sigma^2 + \delta^2)(s/N) - \delta^2 \text{tr}(W), \qquad (4.8)$$

and the first inequality in (4.8) becomes an equality for binary designs. Thus, as our computations also confirm, binary designs play a major role in keeping $\psi_d$ small. From (4.7) and (4.8), a lower bound to the efficiency of *d*, under possible model uncertainty as considered here, is given by

$$\text{eff}_{lb}(\rho) = \frac{(1+\rho)(s/N) - \rho\text{tr}(W)}{\text{tr}(H_d^{-1}) + \rho\{\text{tr}(V_d) - \text{tr}(W)\}}, \qquad (4.9)$$

where $\rho = \delta^2/\sigma^2$, and $V_d$ and *W* are as shown in Lemma 3. Note that *W* does not depend on the design *d*. Just as the bound in (3.7), the one in (4.9) is a conservative benchmark; the actual efficiency of *d* among *N*-run exact designs, in the present setup, is often higher than this bound. The examples in the next section illustrate how procedures A, B1 and B2 introduced earlier can yield designs with impressive values of $\text{eff}_{lb}(\rho)$ even for relatively small run sizes.

We now compare our approach in some detail with that in Wilmut and Zhou (2011), Lin and Zhou (2013) and Yin and Zhou (2014), who explored binary designs under orthogonal parametrization. Instead of assigning a prior on $\xi$, they studied minimaxity under a condition which amounts to

$$\xi^T \xi \leq \delta^2 \qquad (4.10)$$

in our setup. By (4.3) and Lemma 3(a), the counterpart of $\psi_d$, under (4.10), is found to be

$$\tilde{\psi}_d = \max_{\xi : \xi^T \xi \leq \delta^2} \text{tr}\{\text{MSE}_d(\hat{\theta})\} = \sigma^2 \text{tr}(H_d^{-1}) + \delta^2 \lambda_{\max}(V_d - W).$$

In contrast to what happens under an orthogonal parametrization, even with all factors at two levels, the matrix *W* has a complex form in our setup and hence, unlike in (4.7), no further splitting up of $\tilde{\psi}_d$ is possible due to nonlinearity of $\lambda_{\max}$ in its arguments. Moreover, due to lack of differentiability of $\lambda_{\max}$, even the approximate theory does not help in finding a useful lower bound on $\tilde{\psi}_d$. Thus, with condition (4.10), exhaustive or near exhaustive search seems to be the only viable design strategy and, computationally, this may be quite formidable unless $v$ and $N$ are rather small. On the other hand, our condition (4.4) serves essentially the same purpose as (4.10) but, aided by approximate theory, keeps the computations manageable. Incidentally, a connection between $\tilde{\psi}_d$ and $\psi_d$ emerges if we rewrite (4.10) as $\lambda_{\max}(\xi\xi^T) \leq \delta^2$, showing that $\tilde{\psi}_d = \max E_\pi[\text{tr}\{\text{MSE}_d(\hat{\theta})\}]$, where the maximum is over all degenerate priors meeting (4.4). On the other hand, our $\psi_d$ in (4.5) or (4.7) is the corresponding maximum over the wider class all priors, degenerate or not, satisfying (4.4).

## 5. Examples

For ease in presentation, we show the examples in the form of tables. Tables 1-7 exhibit how in a wide variety of situations, approximate theory, in conjunction with procedures A, B1 and B2, yields efficient designs which are also model robust, having $\text{eff}_{lb}$ and $\text{eff}_{lb}(\rho)$ values 0.90 or higher. Since $\text{eff}_{lb}$ and $\text{eff}_{lb}(\rho)$ are only conservative lower bounds, the true efficiencies of these designs, among exact designs with the same run size, should be even better.

Throughout, we take $\rho = 1$ and 5. The choice $\rho = 1$ is along the lines of Yin and Zhou (2014), while the choice $\rho = 5$ sheds light on the robustness of the designs when the departure from the assumed model is possibly even larger. Satisfyingly, our computations show that in binary designs $\text{eff}_{lb}(\rho)$ falls off rather slowly with increase in $\rho$. For non-binary designs, however, the fall is quite fast as anticipated from the findings in the last section. Indeed, all the designs reported in Tables 1-7 are binary.

In each table, we indicate the requirement set $R$ and the value of $q+1$, the smallest run size needed to keep the factorial effect parameters dictated by $R$ estimable under the assumed model. Designs are reported for eight consecutive values of $N$, close to $q+1$. We continue to write $d(N)$ to denote a design with run size $N$. In order to save space, the treatment combinations in a design are listed by their labels following (2.3). Moreover, if one design is similar to another, then we only mention how they differ. Thus in Table 2, the design for $N = 16$ is described as $d(15) + (52, 91) - 4$, to indicate that it can be obtained from $d(15)$ by adding the treatment combinations labeled 52 and 91 and deleting the treatment combination labeled 4.

In each table, after finding the optimal design measure via the algorithm in (3.4) and (3.5), all three procedures A, B1 and B2 are applied for every $N$ considered, and we report the best design also indicating the corresponding procedure. Since none of A, B1 and B2 is very computation intensive, we suggest that this be done in other applications as well. In the tables, whenever procedure A or B1 is mentioned, we indicate the initial designs that they start from.

Table 1. *Robust efficient designs for a $2^6$ factorial*
[$R = \{F_1,...,F_6\} \cup \{F_i F_j : i = 1, 2, 3$ and $j = 4, 5, 6\}$, $q+1 = 16$]

| N | $d(N)$ | Procedure | $\text{eff}_{lb}$ | $\text{eff}_{lb}(1)$ | $\text{eff}_{lb}(5)$ |
|---|---|---|---|---|---|
| 16 | 9, 12, 14, 15, 17, 20, 22, 23, 33, 36, 38, 39, 57, 60, 62, 63 | B2 | 0.9411 | 0.9327 | 0.9256 |
| 17 | $d(16) + 1$ | B2 | 0.9512 | 0.9436 | 0.9371 |
| 18 | $d(17) + 49$ | B2 | 0.9426 | 0.9332 | 0.9250 |
| 19 | $d(18) + 6$ | B2 | 0.9393 | 0.9287 | 0.9194 |
| 20 | $d(19) + 4$ | B2 | 0.9411 | 0.9302 | 0.9204 |
| 21 | $d(20) + 7$ | B2 | 0.9482 | 0.9379 | 0.9285 |
| 22 | $d(21) + 41$ | B2 | 0.9606 | 0.9523 | 0.9444 |
| 23 | $d(22) + 25$ | B2 | 0.9790 | 0.9742 | 0.9695 |

Table 2. *Robust efficient designs for a $2^5 \times 3$ factorial*
[$R = \{F_1,...,F_6, F_1F_6, F_2F_6\}$, $q+1 = 12$; A initiated from $d(304)$ with $\text{eff}_{lb} = 0.9925$; B1 initiated from the full factorial]

| N | $d(N)$ | Procedure | $\text{eff}_{lb}$ | $\text{eff}_{lb}(1)$ | $\text{eff}_{lb}(5)$ |
|---|---|---|---|---|---|
| 13 | 10, 13, 20, 24, 27, 29, 31, 51, 53, 55, 76, 92, 96 | A | 0.9129 | 0.9066 | 0.9019 |
| 14 | $d(13) + (40, 67) - 55$ | A | 0.9303 | 0.9247 | 0.9205 |
| 15 | $d(14) + 4$ | A | 0.9336 | 0.9279 | 0.9235 |
| 16 | $d(15) + (52, 91) - 4$ | A | 0.9460 | 0.9409 | 0.9369 |
| 17 | $d(16) + (8, 17) - 20$ | A | 0.9568 | 0.9524 | 0.9490 |
| 18 | $d(17) + (12, 88) - 76$ | A | 0.9607 | 0.9564 | 0.9530 |
| 19 | 4, 13, 15, 17, 18, 19, 23, 26, 28, 33, 46, 50, 52, 57, 70, 76, 90, 91, 95 | B2 | 0.9604 | 0.9558 | 0.9521 |
| 20 | 2, 4, 6, 7, 8, 19, 21, 34, 36, 37, 41, 44, 49, 60, 70, 71, 74, 79, 88, 90 | B1 | 0.9609 | 0.9561 | 0.9522 |

Table 3. *Robust efficient designs for a $2^2 \times 3^2 \times 4$ factorial*
[$R = \{F_1,...,F_5\}$, $q+1 = 10$; A initiated from $d(280)$ with $\text{eff}_{lb} = 0.9957$; B1 initiated from the full factorial]

| N | $d(N)$ | Procedure | $\text{eff}_{lb}$ | $\text{eff}_{lb}(1)$ | $\text{eff}_{lb}(5)$ |
|---|---|---|---|---|---|
| 14 | 8, 10, 13, 19, 28, 33, 39, 57, 66, 77, 86, 107, 109, 132 | B1 | 0.9300 | 0.9265 | 0.9240 |
| 15 | $d(14) + 137$ | B1 | 0.9443 | 0.9413 | 0.9391 |
| 16 | $d(15) + (9, 25) - 33$ | B1 | 0.9486 | 0.9456 | 0.9434 |
| 17 | $d(16) + 112$ | B1 | 0.9573 | 0.9546 | 0.9526 |
| 18 | 3, 6, 24, 25, 33, 40, 43, 50, 53, 73, 82, 87, 89, 104, 109, 117, 134, 143 | A | 0.9602 | 0.9575 | 0.9556 |
| 19 | $d(18) + 124$ | A | 0.9619 | 0.9592 | 0.9572 |
| 20 | $d(19) + 80$ | A | 0.9626 | 0.9598 | 0.9577 |
| 21 | 8, 9, 13, 19, 25, 28, 38, 39, 47, 57, 66, 72, 77, 82, 86, 99, 107, 109, 112, 128, 137 | B1 | 0.9587 | 0.9554 | 0.9530 |

Table 4. *Robust efficient designs for a $2^8$ factorial*
[ $R = \{F_1,..., F_8, F_1F_2, F_1F_3, F_1F_2F_3\}$, $q+1 = 12$; A initiated from $d(288)$ with $\text{eff}_{lb} = 0.9800$; B1 initiated from the full factorial]

| N | $d(N)$ | Procedure | $\text{eff}_{lb}$ | $\text{eff}_{lb}(1)$ | $\text{eff}_{lb}(5)$ |
|---|---|---|---|---|---|
| 14 | 11, 22, 60, 92, 100, 125, 137, 152, 167, 186, 208, 209, 230, 251 | A | 0.9088 | 0.9063 | 0.9046 |
| 15 | $d(14) + (147, 160) - 152$ | A | 0.9186 | 0.9162 | 0.9145 |
| 16 | $d(15) + (63, 96) - 92$ | A | 0.9267 | 0.9244 | 0.9228 |
| 17 | $d(16) + (58, 89) - 60$ | A | 0.9389 | 0.9369 | 0.9354 |
| 18 | 5, 20, 44, 49, 74, 83, 103, 118, 130, 137, 148, 159, 163, 190, 204, 213, 232, 249 | B1 | 0.9508 | 0.9490 | 0.9477 |
| 19 | $d(18) + (6, 13) - 5$ | B1 | 0.9497 | 0.9477 | 0.9463 |
| 20 | $d(19) + (210, 215) - 213$ | B1 | 0.9530 | 0.9510 | 0.9497 |
| 21 | $d(20) + (134, 191) - 130$ | B1 | 0.9558 | 0.9539 | 0.9525 |

Table 5. *Robust efficient designs for a $3^5$ factorial*
[ $R = \{F_1,..., F_5\}$, $q+1 = 11$; A initiated from $d(306)$ with $\text{eff}_{lb} = 0.9933$; B1 initiated from the full factorial]

| N | $d(N)$ | Procedure | $\text{eff}_{lb}$ | $\text{eff}_{lb}(1)$ | $\text{eff}_{lb}(5)$ |
|---|---|---|---|---|---|
| 15 | 10, 23, 27, 31, 39, 61, 65, 83, 94, 125, 154, 165, 178, 208, 231 | A | 0.9107 | 0.9080 | 0.9062 |
| 16 | $d(15) + (141, 197) - 10$ | A | 0.9203 | 0.9177 | 0.9159 |
| 17 | $d(16) + (114, 239) - 141$ | A | 0.9298 | 0.9273 | 0.9256 |
| 18 | $d(17) + (136, 162) - 154$ | A | 0.9362 | 0.9338 | 0.9321 |
| 19 | 8, 10, 13, 24, 35, 46, 57, 66, 99, 101, 112, 119, 139, 163, 165, 205, 213, 230, 241 | B1 | 0.9378 | 0.9354 | 0.9337 |
| 20 | $d(18) + (7, 46, 66, 100) - (39, 231)$ | A | 0.9429 | 0.9405 | 0.9388 |
| 21 | $d(19) + (74, 106, 167, 207) - (66, 205)$ | B1 | 0.9478 | 0.9455 | 0.9438 |
| 22 | $d(21) + (58, 138) - 139$ | B1 | 0.9512 | 0.9490 | 0.9473 |

Table 6. *Robust efficient designs for a $2^4 \times 3 \times 4$ factorial*
[ $R = \{F_1,..., F_6, F_5F_6\}$, $q+1 = 16$; B1 initiated from $d(74)$ with $\text{eff}_{lb} = 0.9973$, as given by B2]

| N | $d(N)$ | Procedure | $\text{eff}_{lb}$ | $\text{eff}_{lb}(1)$ | $\text{eff}_{lb}(5)$ |
|---|---|---|---|---|---|
| 20 | 1, 24, 28, 29, 50, 63, 73, 93, 105, 123, 133, 134, 160, 161, 162, 163, 164, 166, 167, 180 | B1 | 0.9204 | 0.9156 | 0.9121 |
| 21 | $d(20) + 34$ | B1 | 0.9151 | 0.9098 | 0.9059 |
| 22 | $d(21) + (85, 89) - 29$ | B1 | 0.9167 | 0.9112 | 0.9071 |
| 23 | $d(22) + (101, 157) - 85$ | B1 | 0.9218 | 0.9164 | 0.9123 |
| 24 | $d(23) + (117, 119) - 167$ | B1 | 0.9216 | 0.9159 | 0.9116 |
| 25 | $d(24) + (19, 175) - 163$ | B1 | 0.9250 | 0.9192 | 0.9149 |
| 26 | $d(25) + (54, 138) - 162$ | B1 | 0.9281 | 0.9223 | 0.9179 |
| 27 | $d(26) + (20, 176) - 164$ | B1 | 0.9323 | 0.9266 | 0.9223 |

Table 7. *Robust efficient designs for a $2^4 \times 3^3$ factorial*
[$R = \{F_1,..., F_7, F_1F_2, F_1F_3, F_2F_3, F_1F_2F_3\}$, $q+1 = 15$; B1 initiated from $d(98)$ with eff$_{lb}$ = 0.9959, as given by B2]

| N  | d(N) | Procedure | eff$_{lb}$ | eff$_{lb}(1)$ | eff$_{lb}(5)$ |
|----|------|-----------|-----------|---------------|---------------|
| 20 | 6, 38, 52, 77, 90, 91, 124, 137, 159, 184, 192, 224, 237, 256, 271, 314, 342, 353, 389, 412 | B1 | 0.9202 | 0.9184 | 0.9171 |
| 21 | d(20) + (343, 356) − 353 | B1 | 0.9271 | 0.9253 | 0.9240 |
| 22 | d(21) + (13, 215) − 184 | B1 | 0.9385 | 0.9369 | 0.9358 |
| 23 | d(22) + (1, 175) − 13 | B1 | 0.9416 | 0.9400 | 0.9389 |
| 24 | d(23) + (247, 253) − 256 | B1 | 0.9440 | 0.9424 | 0.9413 |
| 25 | d(24) + (115, 121) − 124 | B1 | 0.9461 | 0.9445 | 0.9433 |
| 26 | d(25) + (49, 206) − 215 | B1 | 0.9485 | 0.9468 | 0.9457 |
| 27 | d(26) + (83, 99) − 90 | B1 | 0.9521 | 0.9505 | 0.9494 |

For smaller $N$ like the ones in the tables, procedure A often yields binary designs and completes well with B1 even under model misspecification. Procedure B2 tends to perform worse though there are exceptions as in Table 1. For $N$ larger than those in the tables, B1 and B2 yield designs with even higher eff$_{lb}$ and eff$_{lb}(\rho)$. Thus, in the setup of Table 5 with $N = 28$, procedure B1, initiated from the full factorial, yields a design with eff$_{lb}$ =0.9608, eff$_{lb}(1) = 0.9584$ and eff$_{lb}(5) = 0.9567$, while in the setup of Table 6 with $N = 33$, procedure B2 leads to a design with eff$_{lb} = 0.9713$, eff$_{lb}(1) = 0.9682$ and eff$_{lb}(5) = 0.9657$. For such larger $N$, however, procedure A often yields non-binary designs which are very efficient under the assumed model but perform poorly under model misspecification.

The smallest $N$, say $N_0$, in Tables 2-7 is a little larger than $q+1$. For $q+1 \leq N \leq N_0 - 1$, in these examples none of A, B1 or B2 yields designs having eff$_{lb}$ and eff$_{lb}(\rho)$ values 0.90 or higher. However, since the bounds are conservative, this does not necessarily mean that these procedures, especially A and B1, lead to inefficient designs for such values of $N$. Simply, the bounds do not assure high efficiencies of these designs and complete enumeration seems to be the only way of assessing their efficiencies. This is again infeasible because $v$ is quite large in these examples. To give a flavor of what may actually happen when $N$ is very close to $q+1$ and the efficiency bounds are not impressive, we discuss below a few situations where $v$ is small and hence complete enumeration is possible.

(a) $2^4$ factorial, $R = \{F_1,..., F_4, F_1F_2, F_3F_4\}$, $q+1 = 7$; $N = 7, 8, 9, 10$;

(b) $2^3 \times 3$ factorial; $R = \{F_1,..., F_4, F_1F_4, F_2F_4\}$, $q+1 = 10$; $N = 10, 11$;

(c) $2 \times 3 \times 4$ factorial; $R = \{F_1, F_2, F_3, F_2F_3\}$, $q+1 = 13$; $N = 13, 14$.

In the situations considered in (a)-(c) above, procedures A, B1 and B2 yield binary designs but none of these has eff$_{lb}$ or eff$_{lb}(\rho)$ values 0.90 or higher. However, a complete enumeration of binary designs shows that A and B1 produce designs with smallest possible $\text{tr}(H_d^{-1})$ as well as smallest possible $\psi_d$ for $\rho = 1$ and 5, in each of these situations except the one in (a) with $N = 9$. In this last situation, both lead to a design with true efficiencies 0.9796, 0.9734 and 0.9664, for $\rho = 0$, 1 and 5, respectively, among all binary designs with the same run size; note that $\rho = 0$ corresponds to the efficiency under assumed model. Thus A and B1 are capable of producing highly efficient, if not optimal, exact designs even when this is not reflected in the values of the corresponding efficiency bounds.

As mentioned earlier, our approach turns to advantage the two features of the baseline parametrization that make the exact design problem hard, namely, nonorthogonality and lack of symmetry among factor levels. Because of these very features, the optimal design measure becomes highly nonuniform and this helps us. Similar approach should be useful also in other situations which exhibit sufficient asymmetry so as to entail such nonuniformity of the optimal design measure. We conclude with the hope that the present work will generate interest in problems pertaining to situations of this kind.

**Appendix**

*Proof of Lemma 1*. By (3.2), the lemma is obvious if either $M(p)$ or $M(\tilde{p})$ is singular. Suppose they are both nonsingular. Then by (2.5) and (3.3), on simplification,

$$M((1-\varepsilon)p + \varepsilon\tilde{p}) = (1-\varepsilon)M(p) + \varepsilon M(\tilde{p}) + \varepsilon(1-\varepsilon)Z^T(\tilde{p}-p)(\tilde{p}-p)^T Z. \quad (A.1)$$

i.e., $M((1-\varepsilon)p + \varepsilon\tilde{p}) - \{(1-\varepsilon)M(p) + \varepsilon M(\tilde{p})\}$ is nnd. Since $M(p)$ and $M(\tilde{p})$ are both nonsingular, the result now follows from (3.2). □

*Proof of Lemma 2*. We proceed along the lines of Silvey (1980, pp. 18-20), with more elaborate arguments to cope with the nonlinearity of $M(p)$ in the elements of $p$. Let $\tilde{p} = (\tilde{p}_1, ..., \tilde{p}_v)^T$ and $p$ be any design measures such that $M(p)$ is nonsingular. The lemma will follow arguing as in Silvey (1980) if we can show that

$$\phi(\tilde{p}) - \phi(p) \geq \lim_{\varepsilon \to 0+} \{\phi((1-\varepsilon)p + \varepsilon\tilde{p}) - \phi(p)\}/\varepsilon. \quad (A.2)$$

and

$$\lim_{\varepsilon \to 0+} \{\phi((1-\varepsilon)p + \varepsilon\tilde{p}) - \phi(p)\}/\varepsilon = \Sigma_{k=1}^v \tilde{p}_k [\text{tr}\{M(p)\}^{-1} - e_k^T\{M(p)\}^{-1}\{M(p)\}^{-1}e_k], \quad (A.3)$$

where $e_k = z_k - Z^T p$, $1 \leq k \leq v$. The truth of (A.2) is not hard to see from Lemma 1. It remains to prove (A.3). To that effect, observe that by (A.1), $M((1-\varepsilon)p + \varepsilon\tilde{p}) = (1-\varepsilon)M(p) + \varepsilon Q(\varepsilon)$, where

$$Q(\varepsilon) = M(\tilde{p}) + (1-\varepsilon)Z^T(\tilde{p}-p)(\tilde{p}-p)^T Z,$$

is nnd, for $0 < \varepsilon < 1$. Write $g = Z^T(\tilde{p}-p)$ and note that by (2.5) and (3.3),

$$M(\widetilde{p}) = \Sigma_{k=1}^{v} \widetilde{p}_k (z_k - Z^T \widetilde{p})(z_k - Z^T \widetilde{p})^T = \Sigma_{k=1}^{v} \widetilde{p}_k (e_k - g)(e_k - g)^T.$$

Therefore, $Q(\varepsilon) = U(\varepsilon)U(\varepsilon)^T$, where $U(\varepsilon)$ consists of the $v+1$ columns $\widetilde{p}_k^{1/2}(e_k - g)$, $1 \leq k \leq v$, and $(1-\varepsilon)^{1/2} g$. Thus for $0 < \varepsilon < 1$, the inverse of $M((1-\varepsilon)p + \varepsilon \widetilde{p})$ is given by

$$(1-\varepsilon)^{-1}[\{M(p)\}^{-1} - \eta\{M(p)\}^{-1}U(\varepsilon)\{I_{v+1} + \eta U(\varepsilon)^T (M(p))^{-1} U(\varepsilon)\}^{-1} U(\varepsilon)^T \{M(p)\}^{-1}],$$

where $\eta = \varepsilon/(1-\varepsilon)$. Hence by (3.2), after some simplification, the left-hand side of (A.3) equals

$$\text{tr}\{M(p)\}^{-1} - \text{tr}[\{M(p)\}^{-1} U_0 U_0^T \{M(p)\}^{-1}],$$

where $U_0 = \lim_{\varepsilon \to 0+} U(\varepsilon)$. Since $\Sigma_{k=1}^{v} \widetilde{p}_k e_k = g$ and, as a consequence, $U_0 U_0^T = \Sigma_{k=1}^{v} \widetilde{p}_k e_k e_k^T$, the truth of (A.3) is now evident. □

**Proof of Lemma 3.** (a) Note the followings facts: (i) As in (3.1), $Z_d^T L_N P_d = Z^T \Delta(r_d) P$. (ii) If $X_d = [1_N \ Z_d]$, then $L_N X_d = [0_N \ L_N Z_d]$, so that as in (3.1), $Z^T \Delta(r_d) X = Z_d^T L_N X_d = [0_q \ H_d]$, using (2.7); therefore, $H_d^{-1} Z^T \Delta(r_d) X = [0_q \ I_q]$. (iii) By the definition of $P$, the matrix $PP^T$ is the orthogonal projector on $C^\perp(X)$, i.e., $PP^T = I_v - X(X^T X)^{-1} X^T$. By (i)-(iii),

$$H_d^{-1} Z_d^T L_N P_d P_d^T L_N Z_d H_d^{-1} = H_d^{-1} Z^T \Delta(r_d) PP^T \Delta(r_d) Z H_d^{-1}$$

$$= H_d^{-1} Z^T \Delta(r_d) \{I_v - X(X^T X)^{-1} X^T\} \Delta(r_d) Z H_d^{-1} = V_d - [0_q \ I_q](X^T X)^{-1}[0_q \ I_q]^T = V_d - W^*,$$

where $W^*$ is the square submatrix of $(X^T X)^{-1}$ as given by its last $q$ rows and columns. Part (a) is now immediate, noting that $W^* = W$, by (2.5) and the definition of $X$.

(b) By (2.5), after some algebra,

$$\Delta(r_d)\Delta(r_d) - \Delta(r_d) = (I_v - N^{-1} r_d 1_v^T)\{D(r_d)D(r_d) - D(r_d)\}(I_v - N^{-1} 1_v r_d^T), \quad (A.4)$$

which is nnd, because $D(r_d)D(r_d) - D(r_d) = \text{diag}(r_{d1}^2 - r_{d1},...,r_{dv}^2 - r_{rv})$ and each $r_{di}$ is a nonnegative integer. Therefore, recalling (3.1), $V_d - H_d^{-1} = H_d^{-1} Z^T \{\Delta(r_d)\Delta(r_d) - \Delta(r_d)\} Z H_d^{-1}$ is nnd. If $d$ is binary, then (A.4) vanishes, and so $V_d = H_d^{-1}$. □

**Acknowledgement.** The work of RM was supported by the J.C. Bose National Fellowship of the Government of India and a grant from the Indian Institute of Management Calcutta. The work of SH was supported by a grant from Kuwait University.